\documentclass[10pt,twoside]{article}
\usepackage{graphicx}
\usepackage{amsmath}
\usepackage{Latex-preprint}
\usepackage{amssymb,latexsym,xy}
\xyoption{matrix}\xyoption{arrow}\xyoption{frame}
\newfam\black
\font\blackboard=msbm10 
\font\blackboards=msbm7
\font\blackboardss=msbm5
\textfont\black=\blackboard
\scriptfont\black=\blackboards
\scriptscriptfont\black=\blackboardss
\def\Bbb#1{{\fam\black\relax#1}}
\def\atitle{}
\def\dual{{\vee}}

\def\CD{{\cal D}}
\def\CI{{\cal I}}
\def\CP{{\cal P}}
\def\CT{{\cal T}}
\def\BC{\Bbb{C}}

\def\BP{\Bbb{P}}
\def\P{\Bbb{P}}
\def\BR{\Bbb{R}}

\def\BZ{\Bbb{Z}}
\def\mapr{\mathop{\longrightarrow}\limits}

\def\TFT{{\rm TFT}}
\def\Coh{{\rm Coh}~}

\def\Mod{{\rm Mod}~}

\def\Sp{{\rm Sp}}

\def\grade{\varphi}

\def\Im{{\rm Im}~}
\def\ch{{\rm ch}}
\def\codim{{\rm codim}}
\def\Hom{{\rm Hom}}
\def\Ext{{\rm Ext}}
\def\Stab{{\rm Stab}}

\def\CF{{\cal F}}

\def\CM{{\cal M}}
\def\CO{{\cal O}}

\def\II{{II}}

\def\CY#1{{CY$_#1$}}
\def\vev#1{{\langle#1\rangle}}

\title{\bf  Dirichlet branes, homological mirror symmetry, 
and stability\vskip -2mm
\vskip 6mm}

\author{Michael R. Douglas\vspace*{-0.5cm}\thanks{
Dept. of Physics and Astronomy,
Rutgers University, Piscataway NJ 08855 USA;
IHES, Bures-sur-Yvette France 91440;
and Isaac Newton Institute, Cambridge CB3 9AL UK.
E-mail: mrd@physics.rutgers.edu}}
\date{June 27,2002}

\begin{document}
\maketitle

\thispagestyle{first} \setcounter{page}{1}

\begin{abstract}\vskip 3mm
We discuss some mathematical conjectures which have come out of the
study of Dirichlet branes in superstring theory, focusing on the case
of supersymmetric branes in Calabi-Yau compactification.  This has led
to the formulation of a notion of stability for objects in a derived
category, contact with Kontsevich's homological mirror symmetry
conjecture, and ``physics proofs'' for many of the subsequent
conjectures based on it, such as the representation of Calabi-Yau
monodromy by autoequivalences of the derived category.

\vskip 4.5mm


\end{abstract}

\vskip 12mm

\section{Introduction} \label{section 1}\setzero
\vskip-5mm \hspace{5mm } 

Let $M$ be a complex manifold with K\"ahler metric,
and $E$ a vector bundle on $M$.  The hermitian Yang-Mills (HYM) equations
are nonlinear partial differential
equations for a connection on $E$, written in
terms of its curvature two-form $F$:
\begin{eqnarray}
F^{2,0} &= 0 \\
\omega \cdot F^{1,1} &= const. \label{eq:hym}
\end{eqnarray}
The theorems of Donaldson \cite{Don:YM} and Uhlenbeck-Yau 
\cite{UY:YM} state that 
irreducible connections solving
these equations are in one-to-one correspondence with $\mu$-stable
holomorphic vector bundles $E$.  

In this talk, we will state a conjecture which generalizes this statement,
motivated by superstring theory, and some results which follow.
In \cite{Doug:ECM} we gave an introduction for
mathematicians to the family of problems which comes out of the study
of Dirichlet branes in type \II\ string compactification on a
Calabi-Yau threefold (henceforth, \CY3), {\it i.e.}  a three complex
dimensional Ricci-flat K\"ahler manifold.  These manifolds are of
central importance in superstring compactification and their study has
led to many results of physical and mathematical importance.

The primary mathematical result about these manifolds, which is the
starting point for their discussion in superstring theory, is Yau's
theorem \cite{Yau}, which states that a complex manifold $M$ with trivial
canonical bundle, or equivalently $c_1(TM)=0$, and which admits a
K\"ahler metric, admits a unique Ricci-flat K\"ahler metric in each
K\"ahler class.  Thus, a family of Ricci-flat metrics is parameterized
locally by a choice of complex structure, and a choice of K\"ahler class.

Given a manifold $M$, one can define a nonlinear sigma model with
target space $M$, using methods of quantum field theory.  It has
been proven to physicists' satisfaction that for $M$ Calabi-Yau,
this sigma model exists and is a $(2,2)$ superconformal field theory (SCFT),
which can be used to compactify superstring theory.  The space of
such theories is parameterized locally by a choice of complex structure
on $M$, and a choice of ``stringy K\"ahler structure.''  A subset of the
stringy K\"ahler structures are specified by a complexified K\"ahler
class, i.e. a point in $H^2(M,\BC)/H^2(M,\BZ)$.  However this is only
a limit (usually called ``large volume limit'') in the full moduli space
$\CM_K$ of stringy K\"ahler structures \cite{AGM,Witten:LSM}; 
we define this space below.

The most well-known results in this area concern mirror symmetry, a
relation between a pair of \CY3's $M$ and $W$
\cite{CK:mbook,Voisin}.  This relation is simple to state in
physics terms: the nonlinear sigma models associated to $M$ and $W$
are the same, up to an automorphism which acts on the $(2,2)$
superconformal algebra.  This implies that every physical observable
for the $M$ model has a mirror observable in the $W$ model, which in a
suitable sense must be the same.  Physics also provides some
techniques for computing these observables, and predicts certain
examples of mirror pairs.

Many interesting and mathematically precise conjectures have been made
based on this equivalence, starting with the work of Candelas {\it et
al} \cite{CDGP} which gave a formula for the number (suitably defined)
of rational curves of given degree in $M$ a quintic hypersurface $M$
in $\BP^4$, in terms of data about the periods of the holomorphic
three-form on its mirror $W$.  This formula has since been proven by
Givental \cite{Givental}, by a tour de force computation 
(see also \cite{LLY}).

The first conjecture which tried to explain mirror symmetry on a
deeper level was made by Kontsevich \cite{Kon:mir}, who proposed that
it should be understood as an equivalence between the derived category
of coherent sheaves on $M$, and the Fukaya category on $W$
\cite{Fukaya}, a category whose objects are isotopy classes of
Lagrangian submanifolds carrying flat connections, and whose morphisms
are elements of Floer cohomology.

It took a little while for physicists to appreciate this conjecture.
For one thing, the natural physical objects to which it relates,
Dirichlet branes, were almost completely unknown when it was made.
D-branes quickly moved to center stage
with the ``second superstring revolution,'' and an open string version
of the Candelas {\it et al} conjecture was made \cite{Vafa},
but the larger physical picture relating to Kontsevich's ideas is
rather subtle and took longer to uncover.

Direct contact was made in \cite{Doug:DC}, where it was shown that the
derived category of coherent sheaves on $M$, $D(\Coh M)$, could be
obtained as a category of boundary conditions in the B-type
topologically twisted sigma model on $M$.  Further developments in
this story appear in \cite{AL:DC,Dia:DC,Laz:DC,AD:DC,Asp:flop}.  This
makes contact with the mathematics of equivalences between derived
categories in various ways, which we will describe.

On realizing how the derived category arises in string theory, one
sees new features which were not present in the standard mathematical
treatments.  The derived category of coherent sheaves is $\BZ$-graded,
but it is clear in string theory that the gradings are naturally
$\BR$-valued, and furthermore depend on stringy K\"ahler data.  This
was mentioned in \cite{Doug:ECM}, and we will describe it in 
more detail here.

The simplest physical question about Dirichlet branes is simply
to know which ones exist as physical branes.  More precisely, while
the general argument leading to the derived category tells us that
every B-type brane corresponds to some object in $D(\Coh M)$,
not all of these objects are actually physical B-type branes, and
we would like to know which ones are.

In the large volume limit, the physical branes are the coherent
sheaves associated to solutions of the HYM equations (with certain
allowed singularities), and this question was answered by the DUY
theorems: the physical branes are the $\mu$-stable sheaves.  To
summarize, one begins with the Hitchin-Kobayashi correspondence,
according to which Yang-Mills connections with $F^{2,0}=0$ are in
one-to-one correspondence with holomorphic bundles with a hermitian
metric.  Such a bundle admits an action of the complexified gauge
group $G_\BC$ (say $GL(N,\BC)$), and one wants to know whether a given
orbit of this group contains a solution of (\ref{eq:hym}).  One can
formally regard the equation (\ref{eq:hym}) as setting an
infinite-dimensional moment map for this action to zero, so the
question is directly analogous to that of defining a symplectic
quotient.  In that context, a $G_\BC$-orbit will contain a solution if
it is stable in the sense of geometric invariant theory (GIT).  The
DUY theorems then state that the appropriate notion of stability is
$\mu$-stability, and prove that a stable orbit indeed contains a
unique solution.  The proof itself requires hard analysis, but at the
end one finds that the intuition from GIT is correct.

This answer has the interesting feature that $\mu$-stability and thus
the set of stable objects depends on the K\"ahler class of $M$,
and can change on walls of codimension $1$.  On
the other hand, it does not depend on the overall scale of the
K\"ahler class, reflecting the scale invariance of the Yang-Mills
equations.

Superstring theory leads to ``generalized HYM equations,''
which reduce in the large volume limit to the HYM equations, but
which explicitly introduce a length scale, the ``string length.''
The next step in the discussion would seem to be to study these
equations.  However, although some results are known in this direction
\cite{MMMS,Leung},
a usable general equation incorporating all stringy effects is not
known at this writing.  We will comment on this problem in
the conclusions.

Rather, what we will do is propose the stringy analog of
$\mu$-stability, a necessary and sufficient condition for an object
in $D(\Coh M)$ to correspond to a physical boundary condition, called
$\Pi$-stability.  This condition will depend on stringy K\"ahler
moduli, and will describe the generalization of wall crossing and
other phenomena in $\mu$-stability to string theory.  

It turns out that the problem of not having a concrete analog of the
HYM equations can be bypassed to some extent by basing the discussion
on SCFT.  In particular, one can prove the ``easy'' half of the DUY
theorems, that a $\Pi$-unstable object cannot be a physical brane,
directly from SCFT.

Despite the direct analogy, the discussion differs dramatically from
that made in geometric invariant theory, because the derived
category is not an abelian category.  Even the most basic elements
of the standard discussion need to be rethought.  In particular, there
is no concept of subobject in the derived category.  We will
overcome these problems, but by paying a price.  Rather than
formulating a condition which can be tested on an object $E$ at a
point in $\CM_K$, we will give a rule which describes the variation of
the entire set of stable objects at a point $p$ in $\CM_K$, call it
$\Stab_p$, to a new set at another point, say $\Stab_{p'}$.  This rule
depends on the path in $\CM_K$ connecting $p$ and $p'$, but
conjecturally only on its homotopy class.

This leads to an interesting phenomenon when one traverses a
homotopically nontrivial loop in moduli space.  Physically, the result
of this must be to recover the original list of stable objects.
However, the specific objects in the derived category which represent
these stable objects can change.  In other words, such a loop is
naturally associated to an autoequivalence in $D(\Coh M)$.
Considerations of mirror symmetry had previously led to conjectures
relating loops and monodromies to autoequivalences
\cite{Kon:mon,ST:braid,Horj:EZ}; 
using $\Pi$-stability we can in a
sense prove these conjectures, as we will explain.

\section{String theory origin of the formalism} \label{section 2}
\setzero\vskip-5mm \hspace{5mm }

We will not try to explain the stringy aspects of the problem in any
detail, but just convey the central ideas.  The most central idea is
the relation between the ``physical'' string world-sheet, described by
SCFT, and the related topologically twisted string or TFT
\cite{Witten:top}.  Given an SCFT, a related TFT is obtained by
choosing a world-sheet supercharge, an operator $Q$ in the $N=2$
algebra satisfying $Q^2=0$, and restricting the Hilbert space to its
cohomology.  In $(2,2)$ SCFT one can make two choices of supercharge,
leading to the ``A-type'' and ``B-type'' TFT's.  In this language, a
mirror pair of \CY3's is a pair $(M,W)$ such that their twisted sigma
models satisfy $\TFT_A(M) = \TFT_B(W)$ and $\TFT_B(M) = \TFT_A(W)$.

The ``stringy'' geometry of the CY $M$ is fully encoded in the SCFT,
while the TFT only gives a small subset of this.  The closed string
sector generally corresponds to metric properties, and the A and B-type
TFT's respectively encode the stringy K\"ahler structure and the complex
structure.  Thus, mirror symmetry can be used to give a definition
to the stringy K\"ahler structure: the stringy K\"ahler moduli space
$\CM_K(M)$ is isomorphic to the standard complex structure moduli space
$\CM_c(W)$.

On general grounds, boundary conditions in an SCFT are objects in a
category, whose morphisms $\Hom(E,F)$ are the $Q$-cohomology of the
Hilbert space of open strings stretching from the boundary condition
$E$ to the boundary condition $F$.  In $(2,2)$ SCFT, this category
will be graded (by ``$U(1)$ charge'').  The category is not guaranteed
to be abelian, but in a sense is generated by a finite set of abelian
subcategories (the branes with specified grading).

A priori, each SCFT (i.e. with given complex and K\"ahler moduli) leads
to a different category of boundary conditions.  A suitable subset of
these descend to TFT, to provide a category of boundary conditions in TFT.
The morphisms are called ``topological open strings.''

In the particular case of the sigma model on $M$, the simplest B-type
boundary conditions are holomorphic submanifolds carrying a
holomorphic bundle \cite{Witten:top}.  If we consider two objects $E$
and $F$ whose support is $M$ and carrying bundles $E$ and $F$, the
operator $Q$ becomes the $\bar\partial$ operator coupled to a
difference of the gauge connections.  Its cohomology, the morphisms
$\Hom(E,F)$ are the Dolbeault cohomology, $H^{0,p}(M,E^\dual\otimes
F)$.  The grading agrees with the standard grading.

More generally, the role of holomorphic bundles in the discussion of
the HYM equations, is played in the string theory generalization by
boundary conditions in $\TFT_B$.  Thus we need to know if this is the
most general boundary condition in $\TFT_B(M)$.  In fact, one can find
direct evidence that it is not, from various sources.  The main point
is the following: the definition of $\TFT_B(M)$ by twisting an SCFT
shows that it is independent of the stringy K\"ahler moduli, and so
one must obtain the same $\TFT_B$ by twisting any SCFT obtained by
deforming in $\CM_K(M)$.  On the other hand, when one compares the
boundary states in examples obtained from various SCFT's, one finds
that they are not the same.

Independence of $\TFT_B(M)$ on $\CM_K(M)$ can only be recovered by
postulating some larger class of boundary conditions which contain the
ones produced by twisting any specific SCFT.  Since the definition of
TFT is made in homological terms, one can attempt to use any of the
standard generalizations known in mathematics.
Kontsevich's conjecture suggests that we try to formulate boundary
conditions which are objects in the derived category of coherent
sheaves on $M$, $D(\Coh M)$.  More generally, we would start with any
of the various abelian categories $A$ of boundary conditions $A$ in
our SCFT, and try to use an object in $D(A)$.  In fact this can be
done, by more or less repeating the standard mathematical construction.

We still need to justify the claim that this is an appropriate
definition of TFT.  This involves two non-trivial steps.  The first of
these is to argue that direct sums of boundary conditions indexed by
an integer $\BZ$ naturally appear as boundary conditions.  One then
sees by more or less standard TFT arguments that the operator $Q$ can
be modified by turning on maps between terms in the complexes, and
that the resulting cohomology is chain maps up to homotopy
equivalence.  The second non-trivial step is to identify complexes
related by quasi-isomorphisms, morphisms which act as the identity on
homology.

The justification of the first step, the $\BZ$-grading on complexes,
is tied to the grading and the structure of the $N=2$ superconformal
algebra \cite{Doug:DC}.  In particular, an A or B-type boundary
condition in SCFT naturally carries a grading.  However, it turns out
that this is not $\BZ$-valued but $\BR$-valued.  This can be used to
make complexes, whose terms have grades differing by integers.  

In fact, it turns out that $D(A)$ is independent of stringy K\"ahler
moduli only as an ungraded category; the gradings depend on the point
in $\CM_K$, and undergo ``flow of gradings.''
We will describe the consequences of this in the next section.

To a certain extent, one can also justify the second step, quotienting
by quasi-isomorphisms, by physical arguments
\cite{Doug:DC,AL:DC,AD:DC}.  The basic intuition is that each brane
$E$ in string theory comes with an ``antibrane'' $\bar E$, an object
whose definition makes clear that a physical configuration containing
a $E$ and a $\bar E$ can ``annihilate'' to a configuration with these
two objects removed.  This can be formalized in the statement that
branes are classified by K theory \cite{Witten:K}, the additive group
completion of the semigroup under direct sum, with deformation
equivalent objects identified.  

One would prefer to use a finer equivalence relation which 
allows determining the holomorphic maps between objects
(the topological open strings).
This requirement can be met by identifying complexes up to
quasi-isomorphism; for example we identify the complex $E_0 \mapr E_1$
with cohomology in one degree, with this cohomology.  Physically, the
intuition behind this is that the brane $E_0$ and the antibrane $\bar
E_1$ can ``form a bound state,'' which is the same physical object
as the cohomology, so this identification is well motivated.
By iterating this construction, one can justify fairly general
identifications \cite{AD:DC}.
Since the derived category is the universal construction of this type
(i.e. makes the minimal identifications meeting the requirements we
just stated), it is clear that one can use it in this context;
we return to this point in the conclusions.

Thus, the boundary conditions in a given SCFT are a subset of the
objects in $D(A)$ for any $A$ obtained from an SCFT related to it
by motion in $\CM_K$, and we seek a condition which specifies
these.  Since the problem reduces to the HYM equations in the large
volume limit, we are led to look for a condition generalizing
$\mu$-stability, which we will describe in the next section.

There is one primary input which can be proven from SCFT and which we
can rely on in formulating this condition.  Namely, we know that our
boundary conditions lead to a unitary $(2,2)$ SCFT, and are governed
by the representation theory of this superconformal algebra.  This leads
to the general result \cite{BFK:CFT} that the grading $\grade$ of a morphism
in a SCFT must satisfy $0 \le \grade \le \hat c = 3$ in our
examples.  On the other hand, flow of gradings would generally lead to
violations of this bound.  This potential contradiction will lead to
the stability condition.

\section{Flow of gradings and $\Pi$-stability} \label{section 3}
\setzero\vskip-5mm \hspace{5mm }

We now try to state the resulting formalism in a mathematically
precise way.  We suppose that we start with a \CY3, $M$, with mirror
$W$.  Recall that $\CM_K(M) \cong \CM_c(W)$, the moduli space of
complex structures on $W$.  The choice of $M$ determines a
``large complex structure limit point'' $p_M\in\CM_K(M)$.  There is
furthermore a ``topological mirror map,'' an isomorphism 
$$
K_0(M) \mapr^m H_3(W,\BZ)
$$
where $K_0(M)$ is the topological K theory group.

A point $p$ in $\CM_c(W)$ determines a holomorphic three-form $\Omega_p$
on $W$, and thus the periods
$$
\Pi: \CM_c(W) \times H_3(W,\BZ) \mapr \BC; \qquad 
 \Pi_p(\Sigma) = \int_\Sigma \Omega_p .
$$
Composing $\Pi \circ m$, we obtain a character on $K_0(M)$ for each $p$,
which we will also call $\Pi_p$.  This is usually called the ``central
charge'' in the physics literature.

The periods $\Pi$ are analytic with singularities on a discriminant
locus $\CD$ of codimension $1$ on which $W$ degenerates.  We define
$\CT_K$ to be the ``stringy Teichm\"uller space,'' a cover of $\CM_K -
\CD$, on which $\Pi$ is single valued.  A deck transformation on
$\CT_K$ corresponds to a linear monodromy on $\Pi$ which preserves the
intersection form, i.e. an element of $\Sp(b_3(W),\BZ)$.

Let $\Stab_p$ be the set of stable objects at a point
$p$ in $\CT_K$, a subset of the objects in $D(\Coh M)$.
Suppose we know $\Stab_p$ at one point $p$; we will
give a rule which, upon following a path $\CP$ in $\CM_K$ from $p$ to
$q$, determines $\Stab_q$.  

We assume that at no point $x$ on the path does $\Pi_x(E)=0$ for an
$E\in\Stab_x$.  Indeed, there are physical arguments \cite{Strominger}
that this can only happen on the discriminant locus (and furthermore
must happen there).

The basic ingredient in the rule is the assignment of $\BR$-valued
gradings to all morphisms between stable objects.  We denote the
grading of a morphism $f\in\Hom(E,F)$ as $\grade_p(f)$.

Given the grading of $f$ between stable objects at a point $p$,
the grading at $q$ obtained by following the path $\CP$ is
\begin{equation} \label{eq:flow}
\grade_q(f) = \grade_p(f) + \Delta_{p,q}\grade(F) - \Delta_{p,q}\grade(E)
\end{equation}
where
\begin{equation} \label{eq:flowobj}
\Delta_{p,q}\grade(E) = \frac{1}{\pi} \Im \log \frac{\Pi_q(E)}{\Pi_p(E)} ,
\end{equation}
resp. for $F$.  We keep track of the branch of the logarithm, so
this defines an $\BR$-valued shift of the grading $\grade(f)$.
We refer to this as ``flow of gradings.''

Given the gradings at an initial point $p$, this determines
the gradings at the endpoint of any path.  Note that flow
preserves the usual additivity of gradings.  Also, morphisms from
an object to itself do not undergo flow.  This implies that Serre
duality is still sensible, and (for our \CY3's) we still have
$\grade(f) + \grade(f^\dual)=3$.

We will assume that the gradings satisfy the condition that if
$\Pi_p(E)/\Pi_p(F)\in\BR$, $\grade(f_{E\mapr F})\in\BZ$.  This is true
in string theory and will be true with the ``initial conditions'' we
describe below.

One can equivalently describe the gradings of morphisms by defining a
grading on stable objects, an $\BR$-valued function $\grade_p(E)$
which satisfies the constraints $\grade_p(E[1]) = \grade_p(E)+1$ and
${\Pi_p(E)}/{e^{i\pi\grade_p(E)}} \in \BR^+$, and varies continuously
on $\CT_K$.  One can then adopt the usual rule that a morphism in
$\Hom(E[m],F[n])$ has grading $n-m$, but now with $m,n\in\BR$.  This
tells us that when an object becomes stable, we only need to specify
the grading of one morphism involving it, to determine the gradings of
all morphisms involving it.

We next consider the distinguished triangles of $D(\Coh M)$.  With
conventional gradings, these take the form
\begin{equation} \label{eq:triangle}
\xymatrix{
&C\ar[dl]^w_{[1]}&\\
A\ar[rr]^u&&B\ar[ul]^v
} 
\end{equation}
and involve two morphisms of degree zero, and
one of degree one, denoted $[1]$.  
After flow of gradings, the gradings of morphisms will typically not
be zero or one, but will always satisfy the relation
\begin{equation}
\grade(u) + \grade(v) + \grade(w) = 1 .
\end{equation}
We now distinguish ``stable'' and ``unstable'' triangles at each
$p\in\CM_K$, according to
their gradings.  In a stable triangle, all three morphisms satisfy the
constraints
\begin{equation}
0 \le \grade(u),\grade(v),\grade(w) \le 1 .
\end{equation}
If any of these constraints is violated, the triangle is unstable.

As one moves along a path in $\CT_K$, the gradings will vary, and the
stability of triangles will change, on walls of real codimension $1$
in $\CT_K$.  Given the assumptions we discussed, one can see that when
the stability of a triangle changes, two morphisms will have grading
zero, and one will have grading one.  This can be seen by plotting the
central charges of the three objects involved in the complex plane.
If we consider the triangle with these vertices, the condition that
the distinguished triangle be stable is that this triangle be embedded
with positive orientation.  A variation of stability is associated to
colinear central charges, and a degeneration of this embedding.

The fundamental rule is that when a stable triangle degenerates, the
object with largest mass $m(E)=|\Pi(E)|$ becomes unstable.
Physically, one says that it ``decays into the other two.''  This is
the string theory generalization of ``wall-crossing'' of
$\mu$-stability under variation of K\"ahler class.  As we discussed
above, this rule follows directly from CFT, independent of any
relation to HYM or other geometric theories.  If we grant the claim
that the HYM equations emerge as the large volume limit of this
formalism, we can say that we have in fact a new derivation of this
phenomenon.

Consistency requires us to implement this rule in the opposite direction
as well: if an unstable triangle with two stable vertices becomes stable,
the third object becomes stable, with grading the same as that of the
other two.  This determines the grading of all other morphisms it
participates in.  For example, when the triangle (\ref{eq:triangle})
becomes stable, the object $B$ becomes stable, with $\grade(B)=\grade(A)$;
this determines the gradings of all morphisms.

Applying this rule in practice requires keeping track of an infinite
set of stable objects and triangles.  Some of the issues involved in
making sense of this are discussed in \cite{AD:DC}.  We conjecture that they
lead to a unique $\Stab_p$ for every $p\in\CT_K$, independent of
choices of order of application of the rules or other ambiguities,
because string theory defines a unique set of B-type branes at each point.

Although we do not have a mathematical proof of consistency, we have
some understanding of it, and can point out two nontrivial ingredients
involved in this claim.  One is the question of whether $\Stab_q$
depends on the choice of path taken from $p$ to $q$.  In general, one
expects such dependence, but only on the homotopy class of path in
$\CM_K-\CD$, which will disappear on $\CT_K$.  On the other hand, one can
construct simple examples which point to a potential problem with this
claim.  Namely, one can find two intersecting walls, along one of
which one has $A$ decaying to $B+C$, and the other with $B$ decaying
to $E+F$.  After crossing the second wall, one is not supposed to
consider the $ABC$ triangle (since $B$ is not stable), leading to the
possibility that following a second path in the same homotopy class
would lead to the contradiction that $A$ remains stable.  This
contradiction does not happen, by an argument given in \cite{AD:DC}
which invokes the ``octahedral axiom'' of triangulated categories
to show that a third triangle $ACF$ will still predict
a decay.  This argument also suggests that there are no ambiguities due
to rule ordering.

The second ingredient is the physical claim that ``the spectrum has
a mass gap,'' in other words that the set
$S_p = \{|\Pi_p(E)|:E\in\Stab_p\}\subset \BR_{\ge 0}$ at fixed $p$
has a positive lower bound.
This assumption prevents infinite chains
of brane decay or creation and is physically widely believed to hold
in these problems.

\section{Results} \label{section 4}
\setzero\vskip-5mm \hspace{5mm }

We first discuss consequences of the physical argument that $D(A)$
should be the same for the various abelian categories $A$ obtained by
explicit SCFT constructions.  First, as one $\CM_K$ can contain
different large volume limits $M$, $M'$, etc., one infers that $D(\Coh
M)\cong D(\Coh M')$ for all of these limits.  In the known physical
examples, the various limits are obtained by taking different moment
maps in symplectic quotients (the linear sigma model construction
\cite{Witten:LSM}), and thus are all birationally equivalent; the
simplest picture is that this equivalence holds for any birationally
equivalent $M$ and $M'$.  Indeed, this is generally believed to be
true (for arbitrary smooth projective varieties) and to some extent
has been proven \cite{Bridgeland,Kawamata}.

One can also define and work with boundary conditions at various other
``exactly solvable'' points in $\CM_K$ such as orbifold points
\cite{DM} and Gepner models \cite{RS}.  These lead to different
abelian categories, some of whose objects can naturally be identified
with holomorphic bundles, and some of which can not.  Physical results on
resolution of orbifolds then lead to the general predictions $D(\Coh
X) \cong D({\rm Coh}^\Gamma \BC^3) \cong D(\Mod-Q_\Gamma)$, where
$\Gamma\subset SU(3)$ is a discrete subgroup, $X\rightarrow
\BC^3/\Gamma$ is a crepant resolution, ${\rm Coh}^\Gamma M$ are
$\Gamma$-equivariant sheaves on $M$, and $Q_\Gamma$ is the path
algebra (with relations) of the McKay quiver for the group $\Gamma$
\cite{Sardo,DGM}.  Indeed this is the subject of the generalized McKay
correspondence \cite{Reid}, and these facts are shown in
\cite{ItoNak,BKR}.

We next turn to discussing $\Pi$-stable objects.
To get started, one needs to know $\Stab_p$ and the gradings at some
point $p$.  The simplest choice is a large volume limit, in which the
stable objects are coherent sheaves.  The gradings which come out of
string theory are not quite the standard gradings of sheaf cohomology,
but are as follows: a morphism $f\in\Ext^n(E,F)$ between sheafs with
support having complex codimension $\codim(E)$ and $\codim(F)$, has
grading $\grade=n+(\codim(F)-\codim(E))/2$.  The central charges $\Pi(E)$ are
in this limit given by
\begin{equation} \label{eq:lvPi}
\Pi(E) = \int_{M} \ch(f_! E) e^{-B-i\omega}\sqrt{\hat A(TM)}
 \qquad + O(e^{-\int\omega}),
\end{equation}
where $\ch$ is the Chern character, $f$ is an embedding
and $f_! E$ is the K-theoretic
Gysin map \cite{MM:K} (one can also evaluate this by using a locally free
resolution and additivity).

Dropping the exponentially small terms (world-sheet instanton
corrections), one can easily show that $\Pi$-stability reduces to a
stability condition for the modified HYM equation considered in
\cite{MMMS,Leung}, and in the limit $\omega >> B,F$ to $\mu$-stability
\cite{DFR:stab}.  Thus $\Pi$-stability does describe B-type branes in
string theory without instanton corrections.

Another limit which is available in certain problems (for example,
orbifolds $\BC^3/\Gamma$ is one in which all central charges $\Pi$ are
real.  One can see that in this case, $\Pi$-stability reduces to
$\theta$-stability in an abelian category \cite{King}, and again this
can be directly justified from string theory \cite{DFR:stab}.  This
simpler description is valid near the limit as long as the
central charges for stable objects live in a wedge $W$ and its
negation $-W$ in the complex plane.

We now give some nontrivial examples in the context of sheaves on
a projective \CY3\ $M$, 
say for definiteness the quintic hypersurface in $\P^4$.
We follow a path in $\CM_K$ parameterized by $J=\int \omega$, decreasing
from large $J$.  In general when $J\sim 1$ the exponential corrections
in (\ref{eq:lvPi})
are important, but their effect can be understood qualitatively
in the examples we consider.

The following is a distinguished triangle in $D(\Coh M)$; 
the numerical superscripts
are the gradings of morphisms at large volume:
$$
\CO \mapr^0 \CO(n) \mapr^{1/2} \CO_\Sigma \mapr^{1/2} \CO[1] \ldots
$$
All three objects are stable at large volume.

Using $\Pi(\CO(n)) \sim \frac{1}{6}(n-iJ)^3$, one finds that 
$\grade(\CO(n)) - \grade(\CO)$ increases as $J$ decreases.  When this
crosses $1$, i.e. $\Pi(\CO(n))/\Pi(\CO)$ is negative real, $\CO_\Sigma$
goes unstable.

If one continues towards the stringy regime, eventually $\Pi(\CO(n))$
and $\Pi(\CO)$ become colinear again.  At this point, the non-sheaf 
$X_n$ defined at large volume by
$$
\CO(n) \mapr^3 \CO \mapr X_n \mapr \ldots
$$
becomes stable.  In \cite{RS}, a Gepner model boundary
state was constructed, which has all the right properties to be this object,
confirming the idea that non-coherent sheaves can appear as stable objects.

Another example uses the sequence
$$
\CI_z \mapr \CO \mapr^{3/2} \CO_z \mapr \CI_p[1] \ldots
$$
where $\CO_z$ is the structure sheaf of a point $z\in M$.
Now the ideal sheaf $\CI_z$ is unstable at large volume, while the other two
objects are stable.  Decreasing $J$ can now decrease the $3/2$,
and eventually the brane $\CI_z$ will become stable.  
This will happen on a line coming into the ``conifold point'' \cite{CDGP}, 
a point near which (choose a local coordinate $\psi$ on $\CM_K$),
$$
\Pi(\CO) \sim \psi.
$$
Continuing around this point further decreases the grade of this map
to $0$, and one finds that the point (or D0) $\CO_z$ becomes unstable.

One can then continue back to large volume.  
Physical consistency requires that a closed
loop $\CP\subset\CM_K$ must induce an autoequivalence $\CF_\CP$, such that
$$
\Stab_{\CP p} \cong \CF_\CP \Stab_p .
$$
Thus we conjecture that $\Pi$-stability leads to such an autoequivalence.

Furthermore, this particular autoequivalence
acts on the K theory class as the mirror Picard-Lefschetz transformation,
$$
[X] \mapr [X'] = [X] + \vev{X,\CO} [\CO] .
$$
These facts and the result we just derived that $\CO_z$ has as
monodromy image $\CI_z$, are enough to prove \cite{AD:DC}
that this autoequivalence
is one conjectured in \cite{Kon:mon,ST:braid,Horj:EZ}, on grounds
of mirror symmetry: $X$ becomes $X'$ defined by the triangle
$$
\CO\otimes \Hom(\CO,X) \mapr X \mapr X' \mapr \dots .
$$
Thus, granting that $\Pi$-stability leads to an autoequivalence,
we have derived the particular autoequivalence associated
to this loop.

Another equivalence of derived
categories which has been obtained from $\Pi$-stability \cite{Asp:flop},
is one shown by Bridgeland to describe a flop \cite{Bridgeland}.  
It is obtained by
going around a point at which the volume of the flopped curve $\sigma$
vanishes,
$$
\Pi_{\sigma} \sim z = B + i J .
$$
This produces the variations of grading 
$$
\Delta\grade(\CO_\sigma) = 1; \ \Delta\grade(\CO_\sigma(-1)) = -1 .
$$
By virtue of
$$
\CO_\sigma(-1) \mapr \CO_\sigma \mapr \CO_z \mapr \ldots,
$$
all points $z\in\sigma$ (and only these points) become unstable
under this variation.  On the other hand, if $\sigma$ is a $(-1,-1)$
curve, one can show that $\dim\Ext^2(\CO,\CO(-1)) = 2$, and the sequence
$$
\CO_\sigma \mapr \CO_\sigma(-1) \mapr X_z \mapr \ldots,
$$
defines ``flopped points,'' which are not sheaves in $D(\Coh M)$, but
can be shown to be precisely the points on the flopped curve $\sigma'$
in the flopped CY $M'$.  Again, this transformation on the
points extends to a unique equivalence of derived categories.

\section{Conclusions and open questions} \label{section 5}
\setzero\vskip-5mm \hspace{5mm }

To summarize, we conjecture that B-type branes in weakly coupled type
\II\ string theory compactified on a \CY3\ $M$, are $\Pi$-stable
objects in $D(\Coh M)$. 

Although the statement of this conjecture relies on string theory, one
can also regard it as {\it defining} a natural class of objects which
are important in mirror symmetry, the stringy generalization of stable
holomorphic bundles (or solutions of the hermitian Yang-Mills
equation), so we think it could be of interest to mathematicians.
Indeed, by assuming the conjecture and simple physical consistency
conditions, one can derive various earlier mathematical conjectures
related to mirror symmetry, as we discussed in section 4.

We believe that $\Pi$-stability as we formulated it in section 3 is
mathematically precise.  In any case the first open question is to
make it precise, prove its consistency, and see to what extent it
uniquely determines the stable objects.  Some results in this
direction have been announced by Bridgeland \cite{Btalk}.

One application, explored in \cite{AL:DC,AD:DC,Asp:flop}, is to
explore ``stringy geometry'' of Calabi-Yau manifolds by using D-branes
as probes \cite{Doug:Dgeom}.  The primary problem of this type would
appear to be the following: at any $p\in\CM_K$, find all ``D0-branes''
(or point objects).  A D0-brane is a family of stable objects $E_z$
with moduli space a \CY3\ $M_p$, such that $D(\Coh M_p)\cong D(\Coh
M)$.  In the second example above, one sees that one can have more
than one D0 ($\CI_z$ and $\CO_z$).  It might be that in other
examples, there is no D0.

A primary application of this conjecture in string theory would be to
determine whether $\Pi$-stable objects with specified $K$ theory class
exist in various regions of $\CM_K$.  In the large volume limit, the
role of $\mu$-stability in this problem was discussed in
\cite{Sharpe}.  For such concrete purposes, it would be better to have
a criterion which can be applied at a single point in $\CM_K$ and does
not require considering the infinite set of stable objects.  There are
examples (e.g. the Gepner point for the quintic) in which it can be
shown that $\Pi$-stability does not reduce to stability in a single
abelian category, but it could be that the stable objects in any small
region could be described as stable objects (in a more conventional
sense) in any of a finite set of abelian categories.

The statement of $\Pi$-stability is admittedly rather complicated, but
this to a large extent only reflects the complexity of its string
theory origins.  It would be interesting to know to what extent it
depends on details of mirror symmetry (such as $M$ a \CY3, and $\Pi$
given by the periods of the mirror $W$), and to what extent it can be
generalized.  The basic definitions make sense for any complex $M$ and
any character $\Pi$ on $K_0(M)$, but we suspect that consistency, and
especially the relation between homotopically nontrivial loops and
autoequivalences, requires further conditions.

Eventually, one would hope to prove the analog of the DUY theorems.
Because general SCFT's on \CY3\ are not rigorously defined, this is
not yet a mathematically well-posed problem.  There are particular
cases such as Gepner models which can be rigorously defined using
vertex operator algebra techniques \cite{Huang}.

Even proving it to the satisfaction of string theorists seems
difficult at present, but might be a reasonable goal.  Such problems
are much studied by string theorists, under the general rubric of
``tachyon condensation.'' \cite{Sen}  The analogy to HYM at least
gives us a guideline for how to proceed in this case.

One would need to start with a working definition of boundary CFT,
including a large enough space of boundary states, and in which
$(2,2)$ supersymmetry is manifest.  One might also start with string
field theory as developed in \cite{Berkovits}.

One could then attempt a direct analog of the Yang-Mills gradient flow
method used in \cite{Don:YM} to construct solutions.  Namely, one
would start with an boundary condition formed as a complex of physical
boundary conditions.  In general this will not be conformal, but will
flow to a conformal boundary condition under renormalization group
flow.  If we start with a large enough class of boundary conditions,
presumably all boundary conditions could be obtained this way.

It seems to us that the first and perhaps the major element in this
project would be to identify and precisely define the group which
plays the role of the complexified gauge group $G_\BC$ in the HYM
theory.  One then needs to address the question of how equivalence
classes of quasi-isomorphic complexes are actually realized in SCFT.
It is not {\it a priori} obvious or necessary for the conjecture to
hold that all such complexes are identified, but only that each stable
class contains a single physical boundary state, which could come
about in various ways.

Finally, one also has the mirror, special Lagrangian description of
Dirichlet branes.  This was also a primary input into formulating the
conjecture (particularly \cite{Joyce}) but this ingredient was to some
extent subsumed by the SCFT considerations.  Nevertheless it is
clearly important to make more connections to this side as well.

We believe that the study of Dirichlet branes on Calabi-Yau manifolds
is proving to be a worthy continuation of the mirror symmetry story,
and will continue to inspire interaction between mathematicians and
physicists for some time to come.

\medskip

I particularly thank my collaborators Paul Aspinwall, Emanuel
Diaconescu, Bartomeu Fiol and Christian R\"omelsberger.  I am also
grateful to Maxim Kontsevich, Greg Moore, Alexander Polishchuk, Paul
Seidel, and Richard Thomas for invaluable discussions.

This work was supported in part by DOE grant DE-FG02-96ER40959.

\end{document}